\documentclass[twocolumn]{article}  
\pdfoutput=1                                                          

\usepackage{amsmath,amsfonts}
\usepackage{graphics,graphicx}
\usepackage{subcaption}
\usepackage{units}
\usepackage{tikz}
\usepackage{ifthen}
\usetikzlibrary{shapes}
\usetikzlibrary{positioning}

\usepackage[english]{babel}

\usepackage[hidelinks]{hyperref}

\usepackage[amsmath]{ntheorem}
\theoremstyle{plain}

\newtheorem{hyp}{Assumption}
\newtheorem{prop}{Proposition}
\newtheorem{lemma}{Lemma}

\newtheorem{rem}{Remark}

\newtheorem{corollary}{Corollary}
\newtheorem{theorem}{Theorem}

\newboolean{LONG} 
\newboolean{CONF} 

\title{\LARGE \bf
Uniform quasi-convex optimisation via Extremum Seeking}

\author{Nicola Mimmo, Lorenzo Marconi and Giuseppe Notarstefano
\thanks{This research was partially supported by the European Project "AerIal RoBotic technologies for professiOnal seaRch aNd rescuE" (AirBorne), Call: H2020, ICT-25-2016/17, Grant Agreement no: 780960.}
\thanks{The authors are with the Department of Electrical and Information Engineering,  "Guglielmo Marconi", University of Bologna, 40126 Bologna, Italy, e-mail: 
        {\tt\small \{nicola.mimmo2, lorenzo.marconi,  giuseppe.notarstefano\}@unibo.it}}%
}

\begin{document}

\setboolean{LONG}{false}
\setboolean{CONF}{true}

\tikzstyle{triangle} = [draw, regular polygon, isosceles triangle]

\maketitle
\thispagestyle{empty}
\pagestyle{empty}

\begin{abstract}
The paper deals with a well-known extremum seeking scheme by proving uniformity properties with respect to the amplitudes of the dither signal and of the cost function. Those properties are then used to show that the scheme guarantees the global minimiser to be semi-global practically stable despite the presence of local saddle points. To achieve these results, we analyse the average system associated with the extremum seeking scheme via arguments based on the Fourier series. 
\end{abstract}

\section{Introduction}

The early researches on the Extremum Seeking (ES) date bake to 1920s \cite{DOCHAIN2011369} and since then this strategy has been extensively exploited to solve several optimisation problems in electronics \cite{Toloue2017Multivariable}, mechatronics \cite{Malek2016Fractional}, mechanics \cite{Zhang2007Numerical}, aerodynamics \cite{Wang2000Experimental}, thermohydraulics \cite{Burns2020Proportional}, and thermoacoustic \cite{Moase2010Newton}.  Some of the most popular ES schemes are those proposed in \cite{tan2006non, KRSTIC2000595}, which represent the subject of the proposed analysis, although a remarkable variety of schemes were proposed, such as the adoption of an integral action in \cite{guay2016perturbation}, the use of a cost function's parameter estimator \cite{Guay2014Minmax,Nesic2013Framework}, the introduction of an observer \cite{HAZELEGER2020109068}, the extension to fractional derivatives in \cite{Malek2016Fractional}, the use of a predictor to compensate output delays in \cite{Oliveira2017Extremum}, the implementation of a Newton-based algorithm avoiding the Hessian matrix inversion \cite{LABAR2019356, Oliveira2020Multivariable}, and the concurrent use of a simplex-method to find the global minimiser \cite{Zhang2016Simplex}.

All the methods, in a way or another, share the common philosophy of perturbing the system  subject to optimisation by means of the so-called \textit{dither} signal, which is periodic in most of the proposed solutions,  to unveil in which direction the associated cost function decreases (in the case of a minimisation problem).  The analysis tool that is typically adopted to prove {\em practical} convergence to the minimiser exploits the \textit{averaging theory}  \cite{Sanders1985Averaging}{. It is shown} that the \textit{average} system {\em asymptotically} converges to the minimiser and that the trajectories of the original system {remain} closed to the \textit{average} ones if a design parameter{,   namely $\gamma$,} is kept sufficiently small  \cite{Teel2003unified,TAN2005550,KRSTIC2000595,TEEL2000317,TEEL1999329}.

For sake of completeness, it is worth mentioning the results in \cite{scheinker2014extremum,ZHU2022109965} in which the value of the cost function is directly used to perturb the phase of the dither signal rather than to estimate the local gradient. 
Lie-derivative arguments, instead of averaging techniques, are then adopted in the analysis.

Within the previous research context, this paper proposes two contributions.

First, {We propose to study the average systems presented in \cite{tan2006non} via Fourier series arguments. This alternative approach allows to handle a class of non-convex cost functions with a unique global minimiser.}
{We show that the average system trajectories converge to a neighbourhood of the minimiser. The size of this neighbourhood is proportional to $\delta$, which is not required to be small.} 
{As} second contribution, {it} is {shown} that the {addition of an} high-pass filter makes $\gamma$ independent {(uniform)} on the value of the cost function. In fact, we {prove} that semiglobal and practical convergence to the minimiser is achieved with the parameter $\gamma$ only depend{ing} on the Lipschitz constant of the cost function in the domain of interest and not on its value. This allows {for} a tuning of the algorithm that preserves good convergence speed even in presence of large domain of attractions. 

	The rest of this paper is organised as follows. Section \ref{sec:Prob} provides the formulation of the problem and reviews the basic ES scheme proposed by \cite{tan2006non}. In this section, we show that the averaging based on the Taylor expansion is not a genuine representative of the ES scheme whose application can be extended to non strictly convex cost functions via a different averaging analysis. Section \ref{sec:Global} describes the ES scheme improved with the high-pass filter and states the result about the uniformity of $\gamma$ with respect to the cost amplitude. The performance of the investigated ES schemes are tested through the simulations detailed in Section \ref{sec:Application}. Finally, Section \ref{sec:Concl} closes this paper with some concluding comments. For lack of space, all the proofs of Lemmas, Propositions, and Theorems claimed in this paper are reported in \cite{Mimmo2022ESarXiv}.

\section{Problem Formulation }
\label{sec:Prob}

The ES problem consists in the optimisation of a{n unknown} cost function $h\,:\,\mathbb{R}\to \mathbb{R}$ {satisfying} the following  two assumptions. 
\begin{hyp}
The function $h$ is smooth and there exists a  $x^\star \in \mathbb{R}$ such that 
	\label{hyp:Existence}
	\[ h(x)-h(x^\star) > 0\quad \forall\, x \in \mathbb{R}\,:\, x \ne x^\star.
	\]
\end{hyp} 
\begin{hyp}	
\label{hyp:Unicity}
		There exist a locally Lipschitz {and \textit{strictly quasi-convex}} function $m\,:\,\mathbb{R}\to \mathbb{R}$, a class-${\cal K}_\infty$ function $\alpha(\cdot)$, and a $A\geq 0$ such that for all $x_1, x_2 \in \mathbb{R}\,:\, (x_1-x^\star)(x_2-x^\star)\ge 0$ \[|m(x_2)-m(x_1)| \ge \alpha(|x_2-x_1|).\]
\end{hyp}
	
	\begin{rem}
 
As for Assumption \ref{hyp:Unicity}, it
asks  a (not necessarily strict) monotone behaviour of $h$, with the latter that could have isolated saddle points. {In other words, $h(\cdot)$ belongs to the class of the so called \textit{strictly quasi-convex} functions \cite{karamardian1967strictly}.} Clearly, Assumption \ref{hyp:Unicity} is weaker than the common assumption $(\partial h(x)/\partial x)x> 0$ for any $x \ne x^\star$ typically present in literature, (see, among the others, [\cite{tan2006non}, Assumptions 3 and 4])  ruling out the existence of local  saddle points. 
	\end{rem}
 
{T}he problem of {\em semiglobal extremum seeking} can be formulated in the following way. For {any  $\epsilon >0$ and $r_0>0$,} design a system of the form 
\[
 \dot x = \varphi_{\epsilon, r_0}(x, h(x),t)  \quad x(0)=x_0, 
\]
so that for all $x_0$ satisfying $|x_0 - x^\star| \leq r_0$ the resulting trajectories $x(t)$ {are bounded and satisfy} $\lim_{t \to \infty}\sup |x(t) - x^\star| \leq \epsilon$.

Among the different ES schemes proposed in literature to solve the previous problem, a common {one} is given by   (see  \cite{tan2006non})
\begin{equation}
	\label{eq:BasicES}
	\dot{x}=  -\gamma \, y_\delta(x,t) \, u(t)\qquad x(0) = x_0
\end{equation}
in which $y_\delta(x,t):=h(x + \delta u(t))$, $u(t):= \sin (2\pi t)$ is the dither signal and $\gamma, \delta > 0$ are tunable parameters\footnote{
The general case of a dither takes the form  $\tau\mapsto\sin(\omega \tau)$ with $\omega >0$, $\tau \in \mathbb{R}$ (as considered in \cite{tan2006non}) can be always obtained by  rescaling the time as $t = \tau \,2\pi/\omega$.}.
 The block diagram of this algorithm is represented in Figure \ref{fig:GloablES}. In the next part we briefly comment the main properties, as available in {the} literature, of this algorithm and strengthen them.  

   Since the right {hand} side of \eqref{eq:BasicES} is 1-periodic, the average system linked to (\ref{eq:BasicES}) is given by (\cite{khalil2002nonlinear}, \S 10.4) 
\begin{equation}
	\label{eq:average}
	\dot{{x}}_a =  -\gamma \int_0^1 y_\delta({x}_a, t) u(t) dt \,.
\end{equation}

In the remaining part of the section we present a different route for the analysis of \eqref{eq:average} based on a Fourier serier expansion rather than on the Taylor one used in the available literature. 
This analysis allows one to claim stronger results on \eqref{eq:BasicES}.

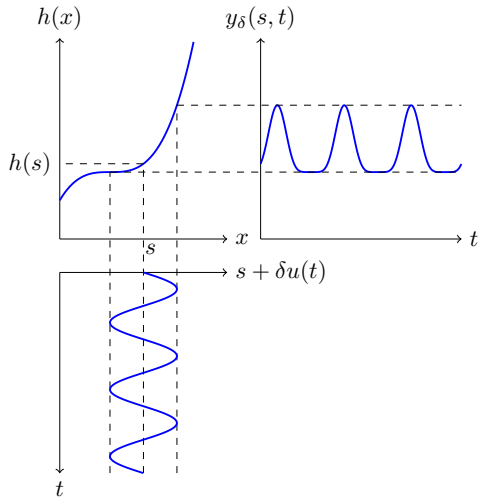
\begin{figure}
	\centering
	\resizebox{0.80\columnwidth}{!}{$
	\begin{tikzpicture}
		\draw[->](0,0) -- (2.5,0) node[right]{$x$};
		\draw[->](0,0) -- (0,3) node[above]{$h(x)$};
		\draw[samples = 100, scale=1, domain=0:2,smooth,variable=\s, blue, thick] plot({\s},{1+(\s-0.75)^3});
		\draw[dashed] ({1.25},{1+(1.25-0.75)^3}) -- (0,{1+(1.25-0.75)^3}) node[left]{$h(s)$};
		\def\xs{3}
		\draw[->](\xs,0) -- (\xs+3,0) node[right]{$t$};
		\draw[->](\xs,0) -- (\xs,3) node[above]{$y_\delta(s,t)$};
		\draw[samples = 100, scale=1, domain=0:3,smooth,variable=\s, blue, thick] plot({\xs + \s},{1+(1.25+0.5*sin(\s*360)-0.75)^3});
		\draw[dashed] ({1.25-0.5},{1+(1.25-0.5-0.75)^3}) -- ({\xs + 3},{1+(1.25-0.5-0.75)^3});
		\draw[dashed] ({1.25+0.5},{1+(1.25+0.5-0.75)^3}) -- ({\xs + 3},{1+(1.25+0.5-0.75)^3});
		\def\ys{-0.5}
		\draw[->](0,\ys) -- (2.5,\ys) node[right]{$s + \delta u(t)$};
		\draw[->](0,\ys) -- (0,-3+\ys) node[below]{$t$};
		\draw[dashed, -] (1.25,-3+\ys) -- (1.25,{1+(1.25-0.75)^3}) node[right, pos = 0.725, xshift = -0.1cm]{$s$};
		\draw[dashed, -] (1.25-0.5,-3+\ys) -- (1.25-0.5,{1+(1.25-0.5-0.75)^3});
		\draw[dashed, -] (1.25+0.5,-3+\ys) -- (1.25+0.5,{1+(1.25+0.5-0.75)^3});
		\draw[samples = 100, scale=1, domain=0:3,smooth,variable=\s, blue, thick] plot({1.25+0.5*sin(\s*360)},{-\s+\ys});
	\end{tikzpicture}
$}
	\caption{Graphical representation of the periodic behaviour of $y_\delta(s,t)$. Conceiving $s$ as a parameter, we see that the oscillation induced by $u(t)$ is elaborated by the non linear map $h(\cdot)$ and results in a periodic function of time (top-right). This latter can be described through its Fourier coefficients $a_k(s)$ and $b_k(s)$.}
	\label{fig:Fourier}
\end{figure}

{The smoothness of $h(\cdot)$ in Assumption 1 is essentially asked to guarantee the existence  of the Fourier series of the function $y_\delta$\footnote{Milder regularity properties guaranteeing the existence of the series could be assumed.}. Then, s}ince $y_\delta(x_a, t)$ and its time derivatives are continuous and periodic, $y_\delta(x_a, t)$ can be expressed in terms of {its} Fourier series  as 
	\begin{equation}
		\label{eq:FourierSum}
		\begin{aligned}
			&y_\delta(x_a, t)= \dfrac{a_{0,\delta}(x_a)}{2} + \\
			&\sum_{k=1}^{\infty}  a_{k,\delta}(x_a) \cos\left(k 2\pi t\right) + b_{k,\delta}(x_a) \sin\left(k 2\pi  t\right)
		\end{aligned}
	\end{equation}
	where
	\begin{equation}
		\label{eq:FourierCoeff}
		\begin{aligned}
			a_{k, \delta}(x_a) &:= 2 \int_{0}^{1} y_\delta(x_a,t)\cos(k 2\pi  t)\,dt\\
			b_{k,\delta}(x_a) &:= 2 \int_{0}^{1} y_\delta(x_a,t)\sin(k 2\pi  t)\,dt\,.
		\end{aligned}
	\end{equation}
     Embedding (\ref{eq:FourierSum}) in (\ref{eq:average}) it is immediately seen that the average system linked to (\ref{eq:BasicES}) reads as 
     \begin{equation} \label{eq:average-Fourier}
      \dot x_a = -{\gamma \over 2} b_{1,\delta}(x_a)\,.
     \end{equation}
  For this system the following result holds. In the result we refer to the class-${\cal K}_\infty$ function $\underline{\delta}^\star(\cdot)$  defined as
  \begin{equation}
  	\label{eq:deltastar}
  \underline{\delta}^\star(s) := 2\int_{0}^{1/2} \alpha(s \sin(2\pi t))\,dt  
  \end{equation}

\begin{lemma} \label{lemma:Reduced1}
Let $h(\cdot)$ be such that Assumptions \ref{hyp:Existence}-\ref{hyp:Unicity} are satisfied. Then:
\begin{itemize}
\item[a)] for all positive $\delta$  there exists a compact set ${{\cal A}_\delta} \subseteq [ x^\star-\delta, x^\star+\delta]$ that is globally asymptotically and locally exponentially stable for (\ref{eq:average-Fourier}).
\item[b)] There exists a $\bar \delta^\star>0$ such that, for all positive $\delta$ such that  $\delta \leq \bar \delta^\star$, there exists an equilibrium point $x_{a \delta}^\star \in \mbox{\em int } {{\cal A}_\delta}$ that is  
 locally exponentially stable for system (\ref{eq:average-Fourier}). If, in addition,  the function $\bar{h}(s) := h(x^\star+s)-h(x^\star)$ is even, then $x_{a \delta}^\star = x^\star$.
\end{itemize}
\end{lemma}
	
Item a) of the previous lemma {states} that the trajectories of the average system reach a compact set ${\cal A}_\delta$ that is contained in a $\delta$ neighbourhood of $x^\star$ for all possible $\delta$. This, in particular, implies that there exists a class $\cal  K L$ function $\beta(\cdot, \cdot): \mathbb{R} \times \mathbb{R}_+ 	\to \mathbb{R}_+$ such that 
\[
|x_a(t)|_{{\cal A}_\delta} \leq \beta(|x_a(0)|_{{\cal A}_\delta}, t)\,
\]
where $|\cdot|_{\mathcal{A}_\delta}$ denote the distance to the set $\mathcal{A}_\delta$. 
Item $b)$ claims that the set ${\cal A}_\delta$ collapses to an equilibrium point if $\delta$ is also taken sufficiently small. Moreover, the last point of item b) shows that $x^\star$ represents the equilibrium point only for  cost functions that are  locally symmetric around the optimum. 

Standard averaging results can be then used to show that {the} same property is preserved also for the trajectories of the original system (\ref{eq:BasicES}) for sufficiently small $\gamma$ but in a semiglobal and practical way. This is detailed in the next Proposition \ref{prop:BasicES} where we refer to the class ${\cal K}_\infty$ function $\chi(s)$ defined as  
\[
\chi(s):= \beta^{-1}(s, 0). 
\]
In the following analysis we  denote by $L_r>0$ and $M_r>0$, {respectively} the local Lipschitz constant and the upper bound  of the function $h(\cdot)$ on a closed interval of length  $r$. In particular, regularity of $h$ implies that   
\begin{itemize}
	\item for all $r>0$ there exists  $L_r > 0$ such that for all $x_1,x_2 \in [x^\star-r,\,x^\star + r]$
	\begin{equation}
		\label{eq:Lr}
		|h(x_1)-h(x_2)|\le L_r |x_1-x_2|\,.
	\end{equation}
	\item for all $r>0$, there exists $M_r > 0$ such that for all $x \in [x^\star-r,\,x^\star + r]$ 
	\begin{equation}
		\label{eq:Mr}
	|h(x)| \leq M_r \,.  
	\end{equation}
\end{itemize}

\begin{prop}
	\label{prop:BasicES}
Let $h(\cdot)$ be such that Assumptions \ref{hyp:Existence}-\ref{hyp:Unicity} hold and let $r,\delta, d$ be arbitrary positive numbers such that $r - d- 2\delta >0$. Let $r_0 :=\chi(r-d-2\delta)$. There exist  a $\bar t(r_0,d)$  and a $\gamma^\star(M_r,L_r,\delta, d) > 0$ such that for any $\gamma \in (0,\,\gamma^\star) $, any $x_0 \in \mathbb{R}\,:\,|x_0-x^\star| \le r_0$,  the trajectories of \eqref{eq:BasicES} are  bounded and
	\[
	 |x(t)|_{{\cal A}_\delta} \le d \qquad \forall \, t \geq {\bar t \over \gamma}\,.
	\]
 \end{prop} 

 An immediate consequence of the previous result is the next corollary showing that under Assumptions \ref{hyp:Existence}-\ref{hyp:Unicity} system (\ref{eq:BasicES}) solves the problem of semiglobal extremum seeking formulated before is fulfilled with a sufficiently small $A$.
\begin{corollary}
	\label{prop:CorollaryBasicES}
Let $h(\cdot)$ be such that Assumptions \ref{hyp:Existence}-\ref{hyp:Unicity} are fulfilled and let $r_0$ and $\epsilon$ be arbitrary positive numbers. Then, there exist $\bar t(r_0, \epsilon)>0$,  $\bar \delta^\star(\epsilon)> 0$ and  $\gamma^\star(r_0,\bar{\delta}^\star, \epsilon) > 0$ such that for any  $\delta \in (0, \bar \delta^\star)$, any $\gamma \in (0,\,\gamma^\star) $ and any $x_0 \in \mathbb{R}\,:\,|x_0-x^\star| \le r_0$,  the trajectories of \eqref{eq:BasicES} are  bounded and
	\[
	 |x(t) - x^\star|  \le \epsilon \qquad \forall \, t \geq {\bar t \over \gamma}\,.
	\]
 \end{corollary} 
 
By going through the proof of Proposition \ref{prop:BasicES},  it is immediately seen that $\gamma^\star$ is  inversely proportional to $M_r$. As a consequence, the higher the cost function is within the set where $x(t)$ ranges,  the lower the value of $\gamma$ and, in turn, the slower the convergence rate of $x$ to the neighbourhood of the optimum. 
 Section \ref{sec:Global} presents an improvement of \eqref{eq:BasicES} overtaking this limitation. 

\subsection{{Comments on Taylor expansion-based averaging analyses}}

The {analysis} of \eqref{eq:average} is typically approached, see \cite{tan2006non}, by using a Taylor expansion of $h(\cdot)$ {to obtain} a system of the form  
\begin{equation}
	\label{eq:av_complete}
	\dot{{x}}_a = \, -\gamma 
	c_1 \delta  \left.\dfrac{\partial h}{\partial x}\right|_{{x}_a} -\gamma 
	\sum_{k=2}^{\infty} c_k \,{\delta^{2k-1}}\,
	\left.\dfrac{\partial^{2k-1}h}{\partial x^{2k-1}} \right|_{{x}_a} 
\end{equation}
where $c_k$ are suitably defined  positive coefficients. A key role in the study of this system is  played by the first order approximation  
\begin{equation}
	\label{eq:First_Order}
	\dot{{x}}_a = \, -\gamma 
	c_1 \delta  \left.\dfrac{\partial h}{\partial x}\right|_{{x}_a}\,. 
\end{equation}
In fact, if {our} Assumption \ref{hyp:Unicity} is strengthen  by asking that{, for any $x \ne x^\star$, it holds} $(\partial h(x)/\partial x)(x-x^\star)> 0$ (respectively $>\alpha(|x-x^\star|)$ with $\alpha(\cdot)$ a class-${\cal K}$ function, see  \cite{tan2006non}, Assumptions 3 and 4){, then the Lyapunov arguments of [\cite{tan2006non}, eq. (45)]} demonstrate $x^\star$ to be a stable (respectively globally asymptotically stable) equilibrium point for (\ref{eq:First_Order}). These stability properties are transferred to \eqref{eq:average} for sufficiently small $\delta$. Then, averaging techniques \cite{khalil2002nonlinear} can be used to prove that{, for sufficiently small $\gamma$,} the trajectories of \eqref{eq:BasicES} {and \eqref{eq:average}} 
 remain arbitrarily close. 
  Namely, the semiglobal extremum seeking problem is solved.   


The fact that the asymptotic properties of the average system \eqref{eq:average} are just ensured by the first order term \eqref{eq:First_Order}, and thus by the gradient of $h$, implies that isolated local saddle points, of $h$, {cannot be handled by that proof technique}. This justifies why [\cite{tan2006non}, Assumption 3], which is stronger than {our} Assumption \ref{hyp:Unicity}, is needed.  Furthermore, we observe that the previous analysis requires that the dither amplitude is kept sufficiently small {for} the higher order terms of the average dynamics {be} negligible.  

\section{{The High-Pass Filer (HPF)}-ES Algorithm}
\label{sec:Global}
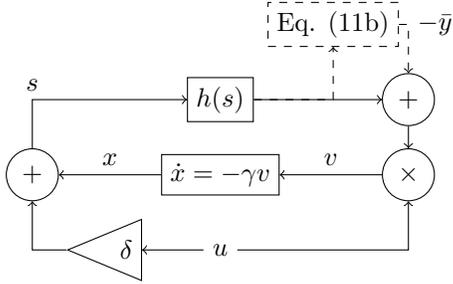
\begin{figure}[t]
	\centering
	\begin{tikzpicture}
		\node[draw](map){$h(s)$};
		\node[draw, right of = map, circle, xshift = 1.5cm](sumy){$+$};
		\node[dashed, draw, above of = map, xshift = 1.5cm](Mean){Eq. \eqref{eq:bary}};
		\draw[dashed, ->] (map) -| (Mean); 
		\draw[dashed, ->] (Mean) -| (sumy) node[pos = 0.5, right]{$-\bar{y}$};
		\node[draw, below of = map, xshift = 0cm, yshift = -0.0cm](ES){$\dot{{x}} = -\gamma v$};
		\node[draw, right of = ES, circle, xshift = 1.5cm](times){$\times$};
		\node[draw, left of = ES, circle, xshift = -1.5cm](sum){$+$};
		\draw[->](map) -- (sumy);
		\draw[->](times) --  node[pos= 0.5, above]{$v$} (ES);  
		\draw[->](ES) -- node[pos= 0.5, above]{${x}$} (sum);
		\draw[->](sum) |- node[pos= 0.5, above]{$s$} (map);
		\node[below of = ES](dither){$u$};
		\draw[->] (sumy) -- (times);
		\node[triangle, left of = dither, shape border rotate=180, xshift = -0.25cm](delta){$\delta$};
		\draw[->](dither) -|  (times);
		\draw[->](dither) -- (delta);
		\draw[->](delta) -| (sum);
	\end{tikzpicture}
	\caption{The HPF-ES scheme consists of an extension of the classic ES algorithm with the addition of the dashed block.}
	\label{fig:GloablES}
\end{figure}

In \cite{tan2006non} the authors proposed the {next} modification of \eqref{eq:BasicES} 
\begin{subequations} 
	\label{eq:GlobalESGeneric}
	\begin{align}
		\label{eq:GlobalESGeneric_dotx}
		\dot x =&\,- \gamma  \left ( y_\delta({x},t)-\bar{y} \right ) \,u(t) &&   {x}(0)={x}_0\\
		\dot{\bar{y}}  = &\, \gamma\left (y_\delta({x},t)-\bar{y}\right ) &&   \bar{y}(0)=\bar{y}_0 \label{eq:bary}
	\end{align}
\end{subequations}
where $u(t)$ is the dither signal defined before, and $(x,\bar{y}) \in \mathbb{R} \times \mathbb{R}$. {A block representation of \eqref{eq:GlobalESGeneric} is depicted in Figure \ref{fig:GloablES}.}

The intuition behind the previous scheme is to interpret $y_\delta({x},t)-\bar y$ as the output of a high pass filter of  $y_\delta$. Moreover, the difference $y_\delta({x},t)-\bar{y}$ represents a proxy  of the local \textit{mean} variation of $h(x)$, directly proportional to the \textit{mean} local Lipschitz constant. In the following we show how this feature guarantees that the upper bound for the value of $\gamma$ is not dependent on $M_r$ {(eq. \eqref{eq:Mr})} but rather only on $L_r$ {(\eqref{eq:Lr})}. As for $\gamma$, similarly to \eqref{eq:BasicES}, it must be small to let $x$ and $\bar{y}$ be sufficiently slow to preserve the correlation between the oscillations of $y_\delta({x},t)-\bar{y}$ and those of $u(t)$.
The average system of \eqref{eq:GlobalESGeneric} is defined as 
\begin{subequations}
	\label{eq:TildeAV}
	\begin{align}
		\label{eq:DotTildeXAv}
		\dot{{x}}_a &=-\gamma  \int_{0}^{1} \left ( y_\delta({x_a},\tau)-\bar{y}_a \right ) \, u(\tau)\,d\tau\\
		\label{eq:baryAV}
		\dot{\bar{y}}_a & = -\gamma\,\bar{y}_a +  \gamma\int_{0}^{1}y_\delta ({x}_a,\tau)\,d\tau.
	\end{align}
\end{subequations}

By expanding  $y_\delta(x_a, t)$ in terms of  Fourier series as in the previous section, it turns out that   
\begin{subequations} 
	\label{rho}
	\begin{align}
		\displaystyle \int_{0}^{1}y_\delta(x_a,\tau)\,d\tau &= \dfrac{a_{0,\delta}(x_a)}{2}\\
		\displaystyle \int_{0}^{1}\left ( y_\delta({x_a},\tau)-\bar{y}_a \right )  u(\tau)\,d\tau  &= \dfrac{b_{1,\delta}(x_a)}{2}
	\end{align}
	where in the latter we exploited $\int_{0}^{1}\bar{y}_a u(\tau) d\tau = 0$ and \eqref{eq:FourierCoeff}. 
\end{subequations}

Hence, the average system reads as 
\begin{subequations}
	\label{eq:TildeAV_2}
		\begin{align}
			\label{eq:DotTildeXAv_2} 
		\dot{{x}}_a &=-\gamma\dfrac{ b_{1,\delta}(x_a)}{2 }\\
		\label{eq:DotBarYAV}		
		\dot{\bar{y}}_a & = -\gamma\, \bar{y}_a + \gamma\, \dfrac{a_{0, \delta}({x}_a)}{2}
	\end{align}
\end{subequations}
which is a cascade where the first subsystem coincides with \eqref{eq:average-Fourier} and the second subsystem is linear and asymptotically stable. From this, the next result follows from  Lemma  \ref{lemma:Reduced1}.

\begin{lemma} \label{lemma:Reduced2}
	Let $h(\cdot)$ be such that Assumptions \ref{hyp:Existence}-\ref{hyp:Unicity} are satisfied. Then:
	\begin{itemize}
		\item[a)] for {any} positive $\delta> 0$ there exist a compact set ${\cal A}_\delta \subseteq [ x^\star-\delta, x^\star+\delta]$ and a continuous function $\tau: \mathbb{R} 	\to \mathbb{R}$ such that the set
		\[
		 \mbox{\rm graph} \left. \tau \right |_{{\cal A}_\delta} = \{ (x_a, y_a) \in {\cal A}_ \delta \times \mathbb{R} \; : \; y_a = \tau(x_a)\}
		\]
 is globally asymptotically and locally exponentially stable for \eqref{eq:TildeAV_2}.
		\item[b)] There exists $\bar \delta^\star>0$ such that, for {any} $\delta \in (0, \bar \delta^\star )$, there exists an equilibrium point $(x_{a \delta}^\star, \bar y_{a \delta}^\star ) \in  \mbox{\rm graph} \left. \tau \right |_{{\cal A}_\delta}$ that is  globally asymptotically and locally exponentially stable for system (\ref{eq:average-Fourier}). If, in addition,  the function $\bar{h}(s) := h(x^\star+s)-h(x^\star)$ is even, then $x_{a \delta}^\star = x^\star$.
	\end{itemize}
\end{lemma}
 
From this, the following result  mimics the one of Proposition \ref{prop:BasicES} with the remarkable difference that the upper bound $\gamma^\star$ of $\gamma$ is uniform with respect to $M_r$.

\begin{theorem}
	\label{prop:HPFES}
	Let $h(\cdot)$ be such that Assumptions \ref{hyp:Existence}-\ref{hyp:Unicity} are fulfilled and let $r,\delta, d$ be arbitrary positive numbers such that $r - d- 2\delta >0$. Let $r_0 :=\chi(r-d-2\delta)$. There exist  a $\bar t(r_0,d)$ and a $\gamma^\star(L_r,\delta, d) > 0$  such that for any $\gamma \in (0,\,\gamma^\star) $ and any $(x_{0},\,\bar{y}_{0})$ fulfilling
	$|x_{0}-x^\star|\le r_0$ and $|\bar{y}_{0} - a_{0,\delta}(x_{0})/2| \le \gamma^\star$,
then	the trajectories of \eqref{eq:GlobalESGeneric} are  bounded and
	\[
	\|(x(t), \bar y(t))\|_{\mbox{\rm graph} \left. \tau \right |_{{\cal A}_\delta} } \le d \qquad \forall \, t \geq {\bar t \over \gamma}\,.
	\]
\end{theorem} 

\section{Numerical Results}
\label{sec:Application}
\begin{figure}
	\centering
	\includegraphics[width=0.55\columnwidth]{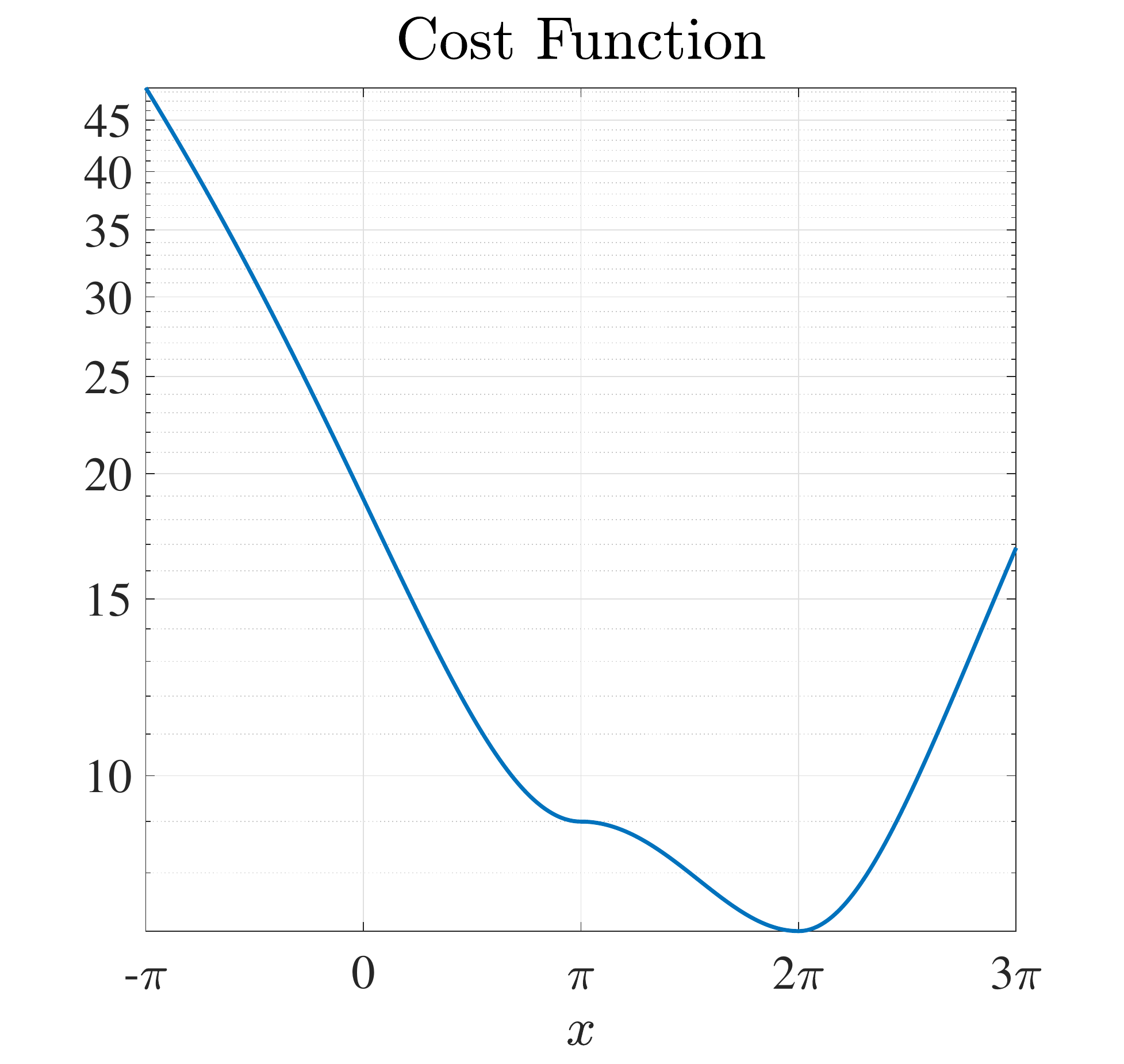}
	\caption{Cost function \eqref{eq:h2ndCase} for $h_0 = 0$. It presents multiple points with null first derivative and it is also non-symmetric around the minimiser. The saddle point is located at $x = \pi$ and the minimiser is at $x^\star = 2\pi$.}
	\label{fig:2ndMap_a}
\end{figure}

This section presents numerical results obtained adopting 

 the following cost function $h(\cdot)\,:\,\mathbb{R}\to \mathbb{R}$ 
\begin{equation}
	\label{eq:h2ndCase}
	h(x) = h_0+\left\{\begin{array}{cc}
		(x-\pi)^2-1 & x < \pi\\
		\cos(x-\pi)-2 & x \in [\pi,\,2\pi)\\
		(x-2\pi)^2-3 & x \ge 2\pi
	\end{array}\right.
\end{equation}
where $h_0 \in \mathbb{R}$. This function, which verifies Assumptions \ref{hyp:Existence} and \ref{hyp:Unicity}, is depicted in Figure \ref{fig:2ndMap_a}, for $h_0 = 10$, {on $[-\pi,\,3\pi]$}. 

\begin{figure}
	\centering
	\includegraphics[width=0.55\columnwidth]{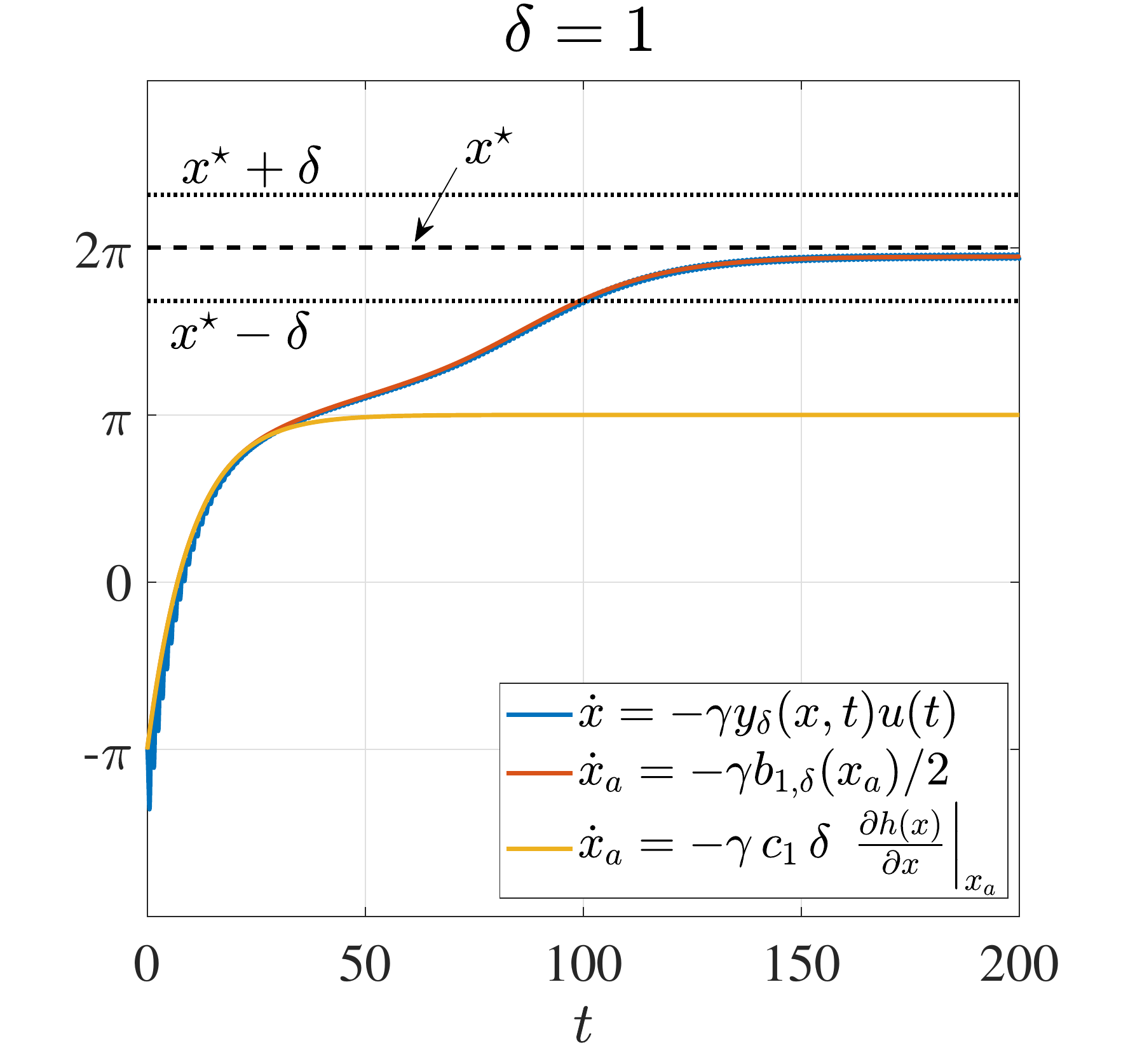}
	\caption{The presence of saddle points stacks the Taylor-based averaging (yellow) of the classic ES (blue). Vice versa, the Fourier-based averaging (red) better represents the actual behaviour of the classic ES (blue). These results are obtained for $h_0=0$ and $\gamma = 0.1$.}
	\label{fig:2ndMap_b}
\end{figure}
\begin{figure}
	\centering
	\includegraphics[width=0.55\columnwidth]{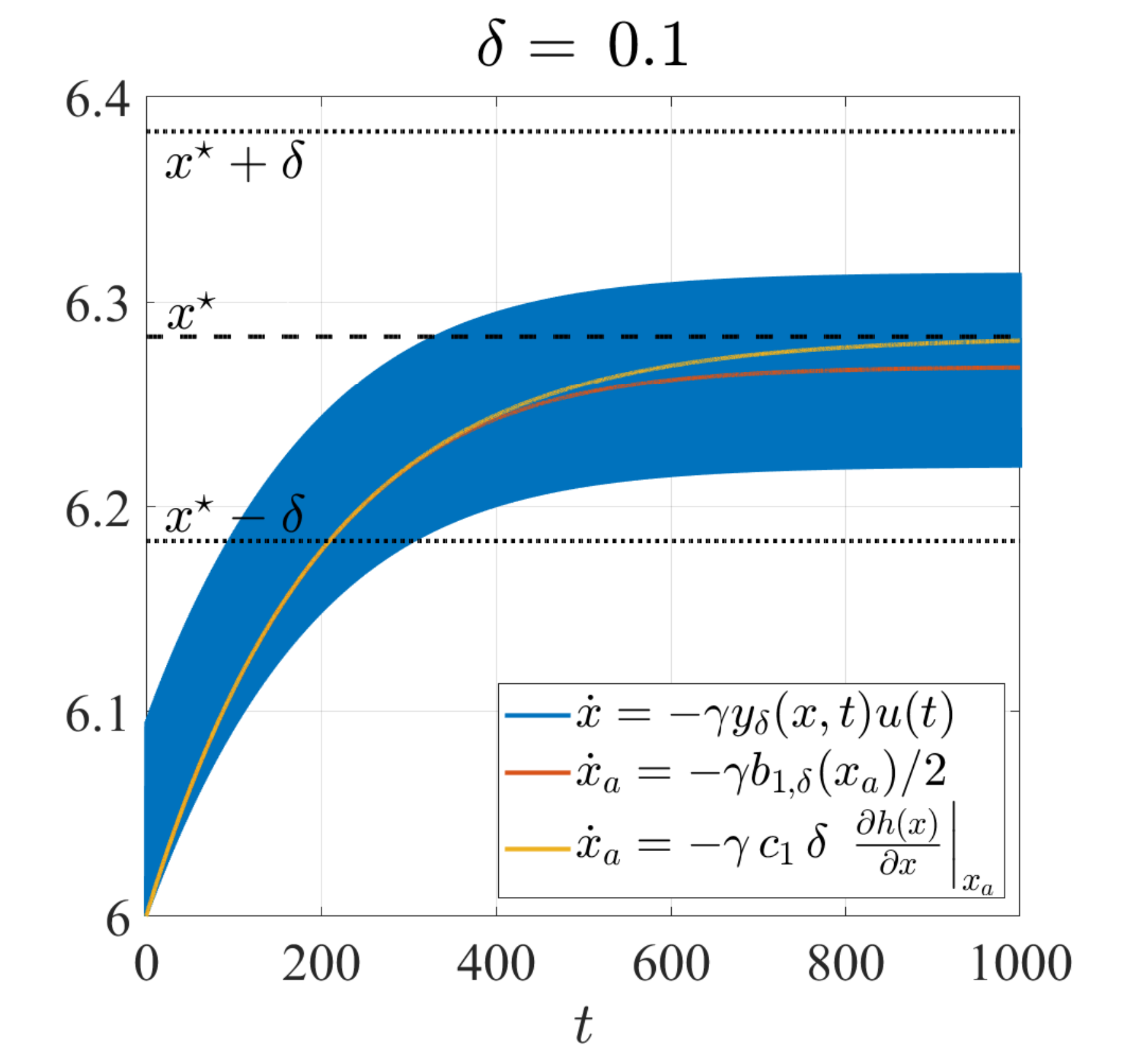}
	\caption{The average based on the Taylor expansion (yellow) of the classic ES wrongly assess $x^\star$ as an equilibrium point. Indeed, due to the asymmetry of $h(\cdot)$ around $x^\star$, the classic ES (blue) converges to 
	an equilibrium point different from $x^\star$. Moreover, the average based on the Fourier series (red) tracks the classic ES more accurately than the average based on the Taylor expansion. These results are obtained for $h_0=0$ and $\gamma = 0.1$.}
	\label{fig:2ndMap_c}
\end{figure}
\begin{figure}
	\centering
	\includegraphics[width=0.55\columnwidth]{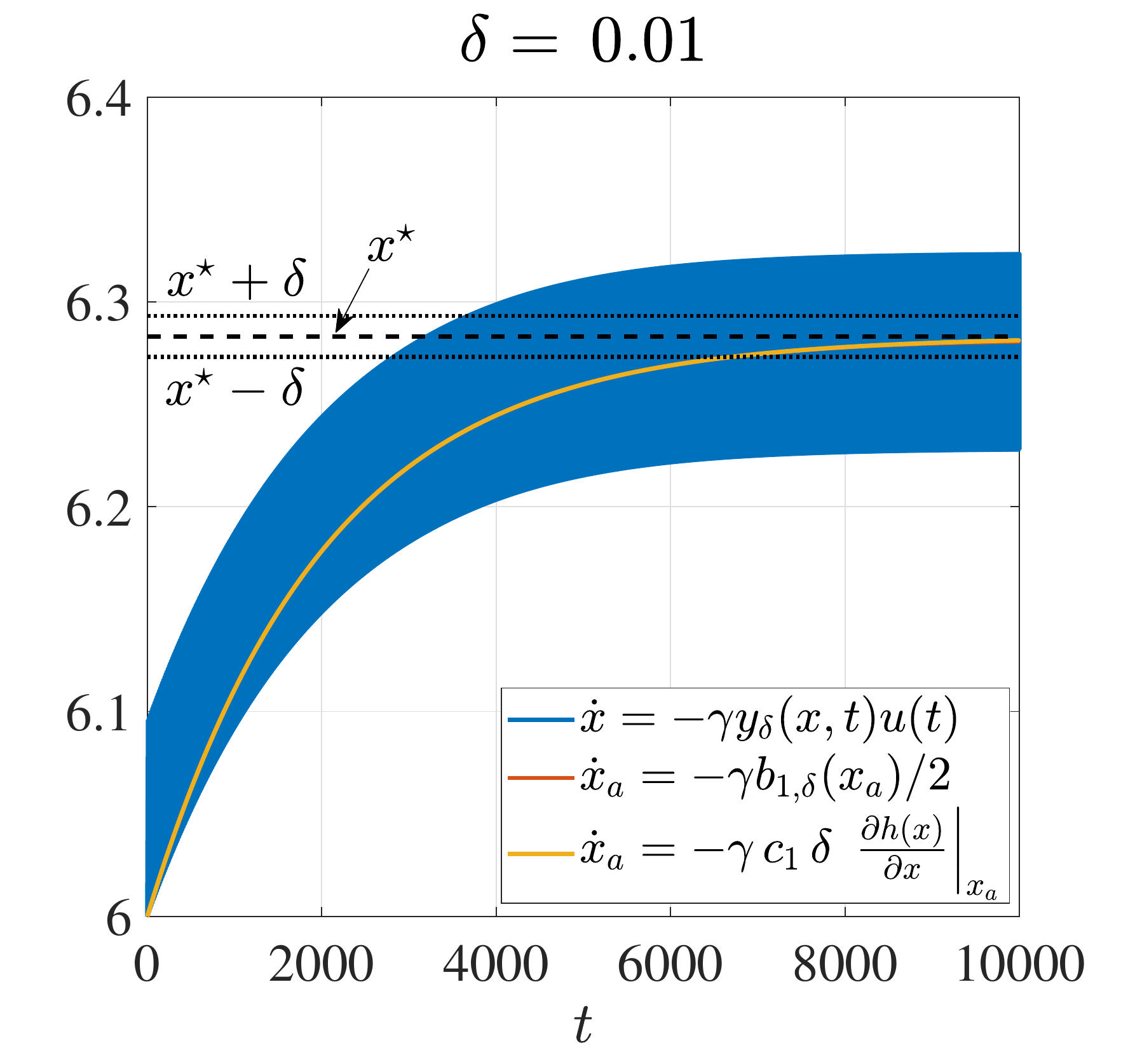}
	\caption{The average of the classic ES, based on the Fourier series (red), converges to an equilibrium point within the set $[x^\star-\delta,\, x^\star+\delta]$ accordingly to what foreseen in Lemma \ref{lemma:Reduced1}. Thus, the smaller is $\delta$ the closer the equilibrium point is to $x^\star$. These results are obtained for $h_0=0$ and $\gamma = 0.1$.}
	\label{fig:2ndMap_d}
\end{figure}

 {Figures \ref{fig:2ndMap_b} and \ref{fig:2ndMap_c} show} that the asymmetry of $h(\cdot)$ in the neighbourhood of $x^\star$ implies that the equilibrium point of \eqref{eq:average} does not correspond to $x^\star$ as wrongly assessed by the averaging through Taylor expansion, see yellow lines in Figures \ref{fig:2ndMap_c} and \ref{fig:2ndMap_d}. {Vice versa, {Figures \ref{fig:2ndMap_c} and \ref{fig:2ndMap_d} show that} the trajectory of \eqref{eq:average-Fourier} (red) tracks, better than those of \eqref{eq:First_Order} (yellow), the trajectory of \eqref{eq:BasicES} (blue) and, as claimed in Lemma \ref{lemma:Reduced1}, converges into a $2\delta$-wide set centred at $x^\star$}.

\begin{figure}[t!]
	\centering
		\begin{subfigure}{.48\columnwidth}
		\includegraphics[width=\textwidth]{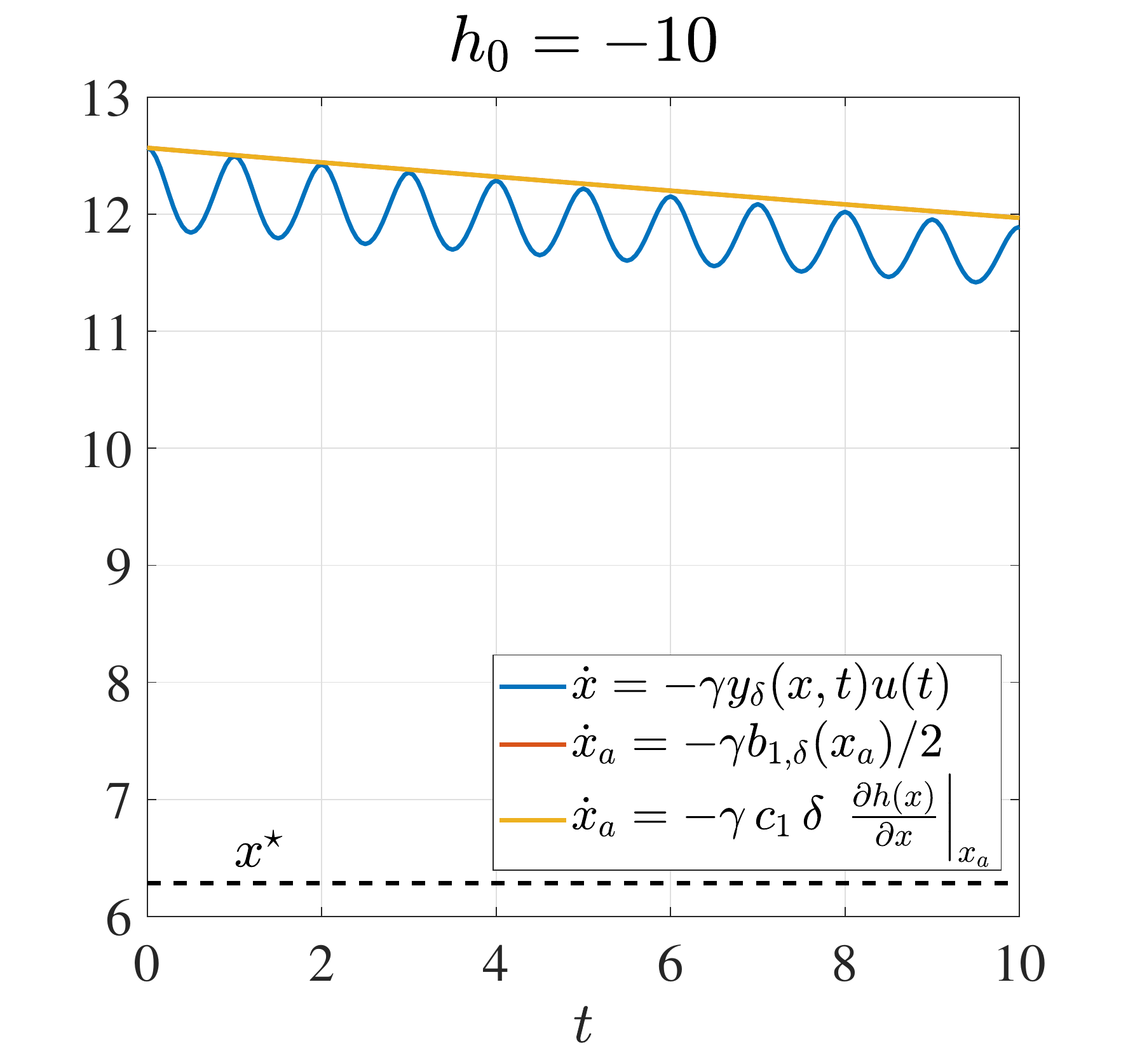}
		\caption{} 
	\end{subfigure}
		\begin{subfigure}{.48\columnwidth}
	\includegraphics[width=\textwidth]{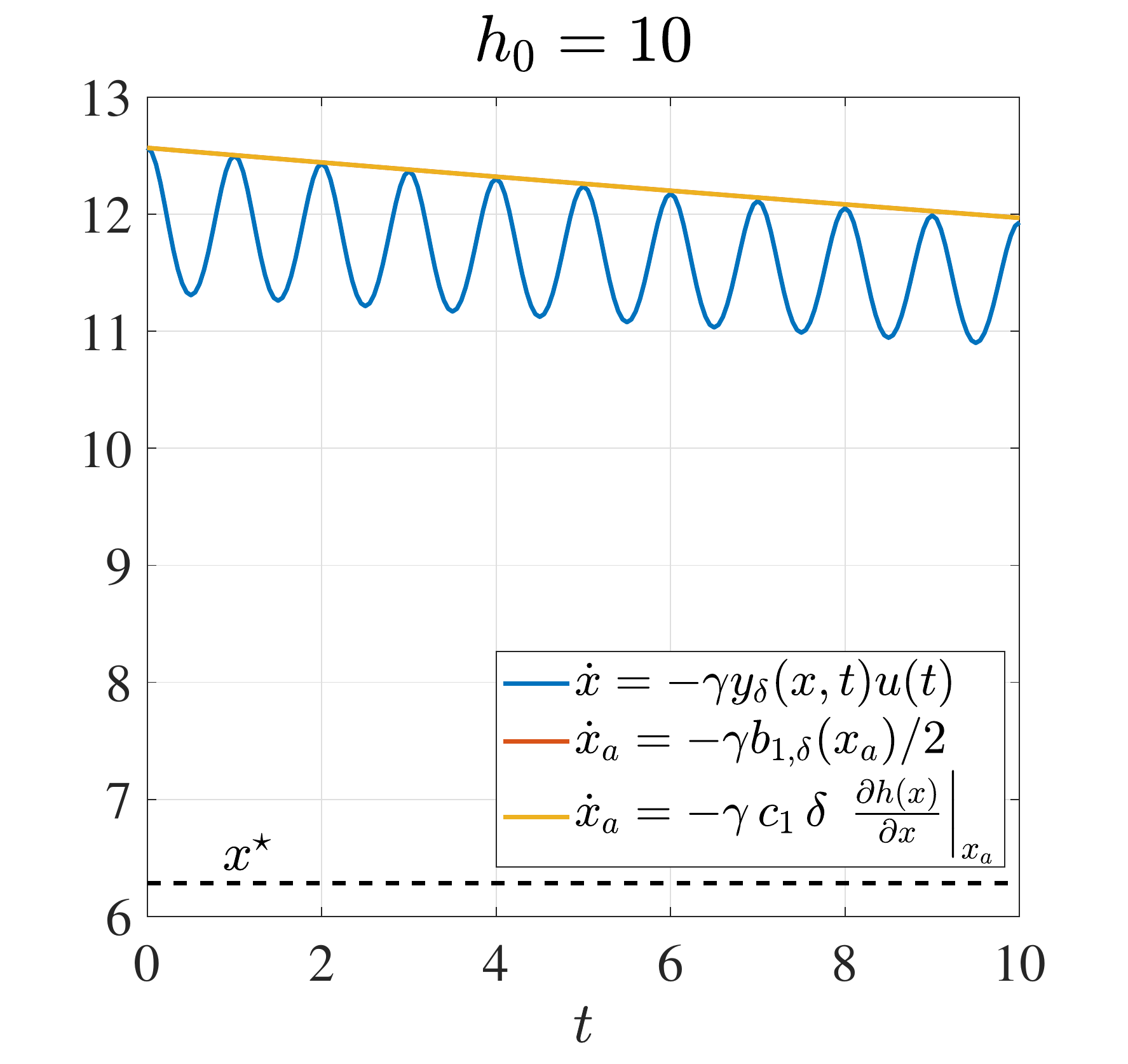}
	\caption{} 
\end{subfigure}
		\begin{subfigure}{.48\columnwidth}
	\includegraphics[width=\textwidth]{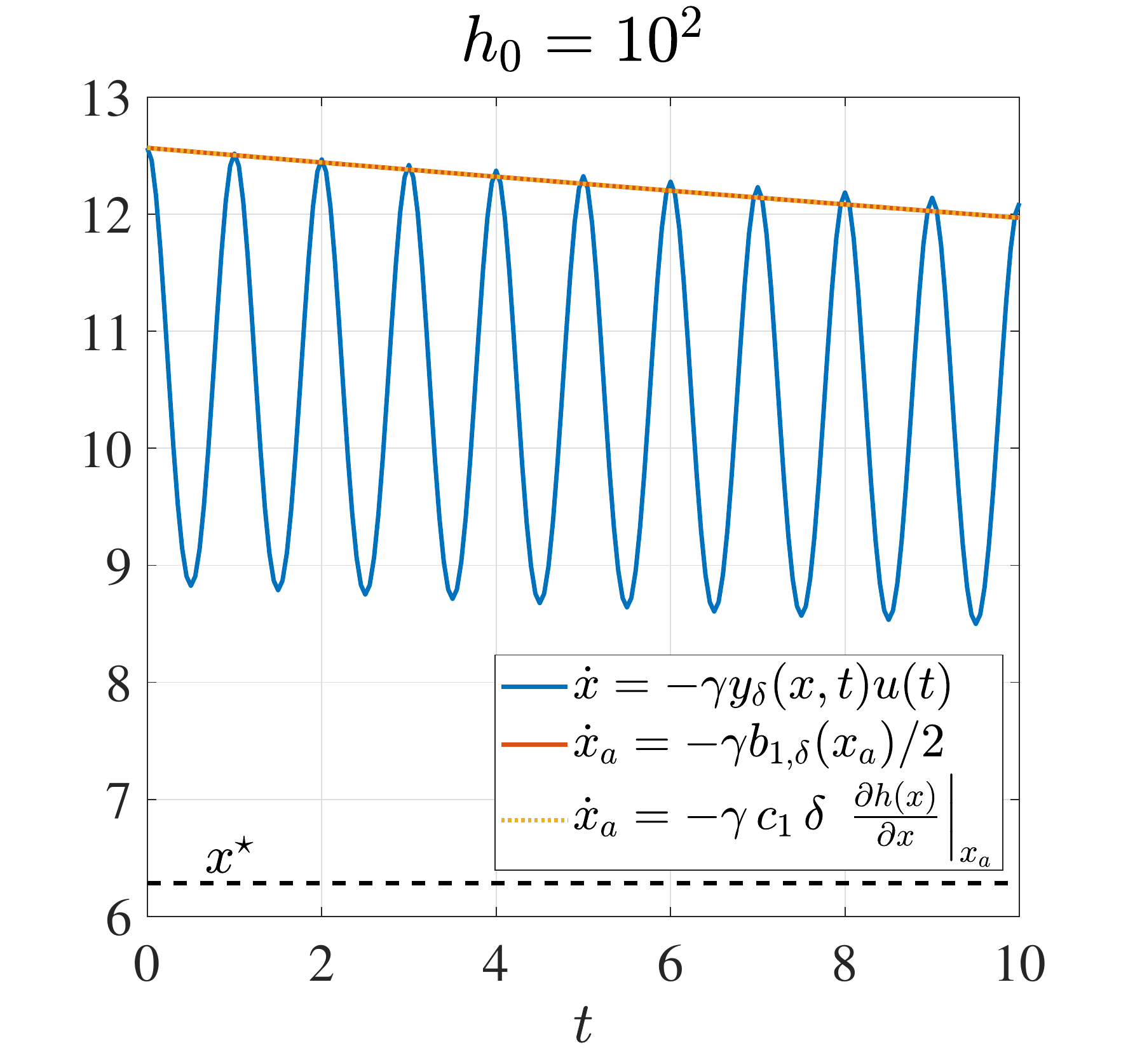}
	\caption{} 
\end{subfigure}
		\begin{subfigure}{.48\columnwidth}
	\includegraphics[width=\textwidth]{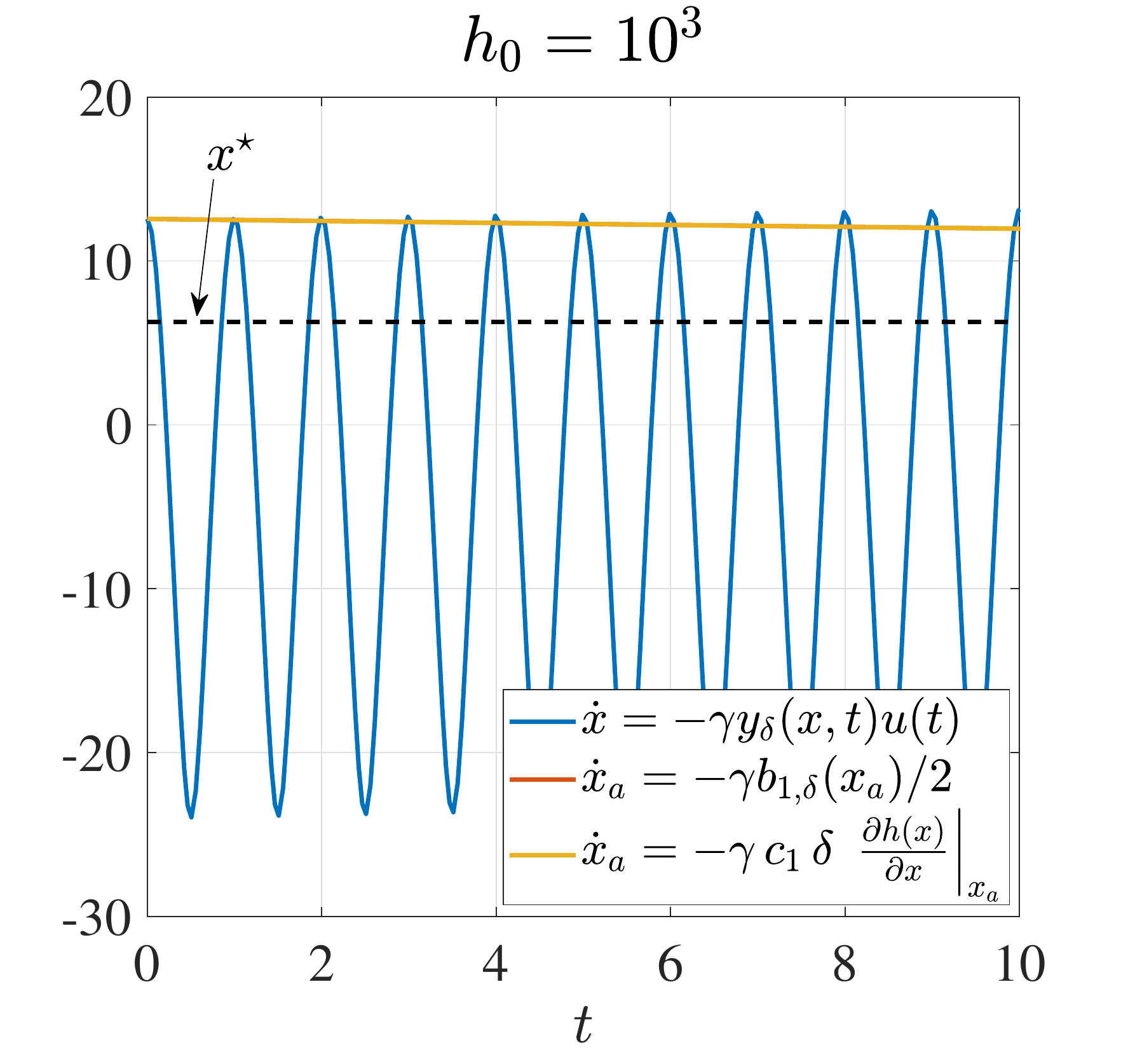}
	\caption{} 
\end{subfigure}
	\caption{The behaviour of the classic ES changes for increasing $M_r$. In these simulations the value of $h_0$ is increased from $-10$ to $10^3$ (subplots from (a) to (d)) demonstrating that, keeping fixed $\gamma$, the classic ES (blue) becomes oscillatory for large values of the cost function. These simulations are performed with $\gamma = \delta = 0.1$.}
	\label{fig:Mr}
\end{figure}

The second test is performed to show that the classic ES suffers of large values of $|h(\cdot)|$. Indeed, as depicted in Figure \ref{fig:Mr} {(blue lines in subplots from (a) to (d))}, while keeping $\gamma$ fixed, larger $M_r$ lead to more oscillatory behaviours. To mathematically support this result we observe that 
\begin{equation}
	\label{eq:Taylor}
	h(x+\delta u) = h(x) + R(x, \delta u)
\end{equation}
where $R(\cdot,\cdot)$ represents the remainder of the Taylor expansion around $x$ of $h(\cdot)$. Exploit the definition of the Lipschitz constant of $h(\cdot)$ and $|u(t)|_\infty \le 1$ to bound the remainder from above as
\begin{equation}
	\label{eq:BoundRemainder}
	|h(x+\delta u) - h(x)| = |R(x, \delta u)| \le L_r\delta.
\end{equation}
Substitute \eqref{eq:Taylor} into \eqref{eq:BasicES} 
\begin{equation}
	\label{eq:ClassicESTaylor}
	\dot{x} = -\gamma h(x+\delta u) u = -\gamma (h(x)+R(x,\delta u)) u
\end{equation}
and investigate the following support system
\begin{equation}
	\label{eq:ClassicESTaylorApprox}
	\dot{x}_1 =-\gamma h(x_1)  u\qquad x(0) = x_{10}
\end{equation}
conceivable as approximation of \eqref{eq:ClassicESTaylor} for $|h(x)| \gg L_r\delta$. Let $H(x):= \int h(x)^{-1} dx$ and solve \eqref{eq:ClassicESTaylorApprox} by parts as
\begin{equation}
	x_1(t) = H^{-1}(H(x_{10}) -\gamma (\cos(t)-1) )
\end{equation}
which is a pure oscillation whose amplitude is proportional to $|H(x_{10})|${, evident in subplot (d) of Figure \ref{fig:Mr}.}

\begin{figure}[h!]
	\centering
	\includegraphics[width=0.48\columnwidth]{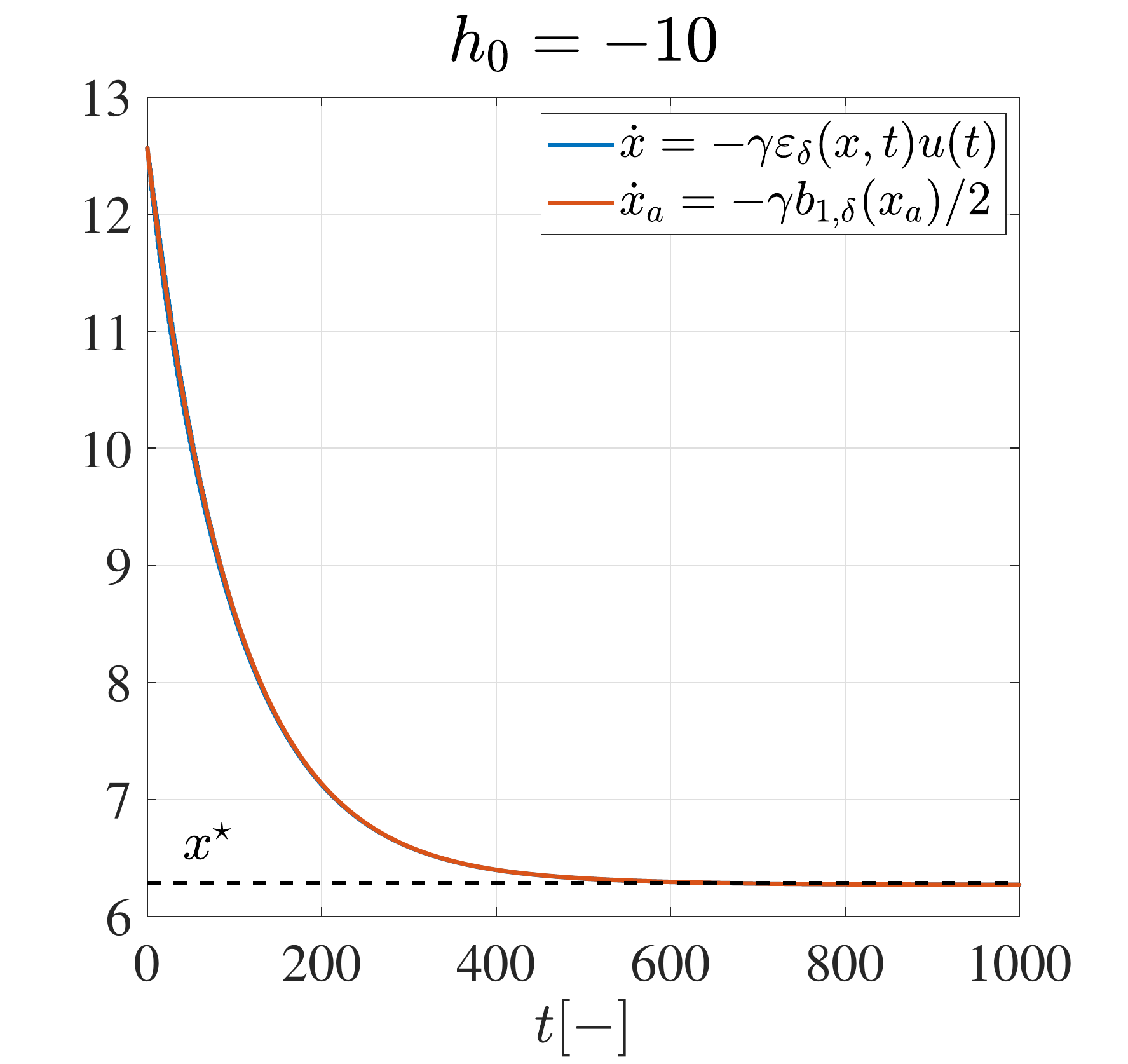}
	\includegraphics[width=0.48\columnwidth]{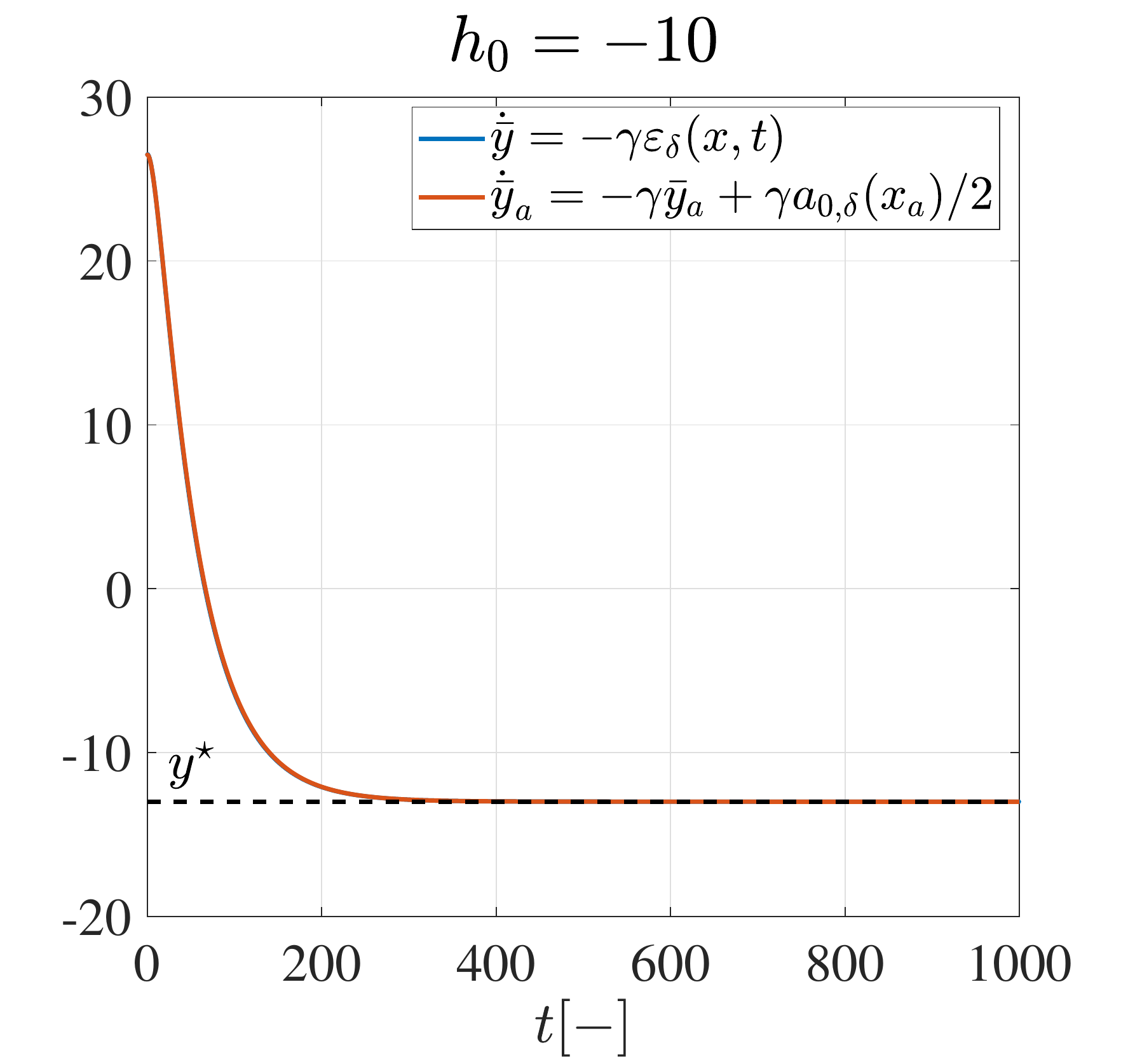}
	\includegraphics[width=0.48\columnwidth]{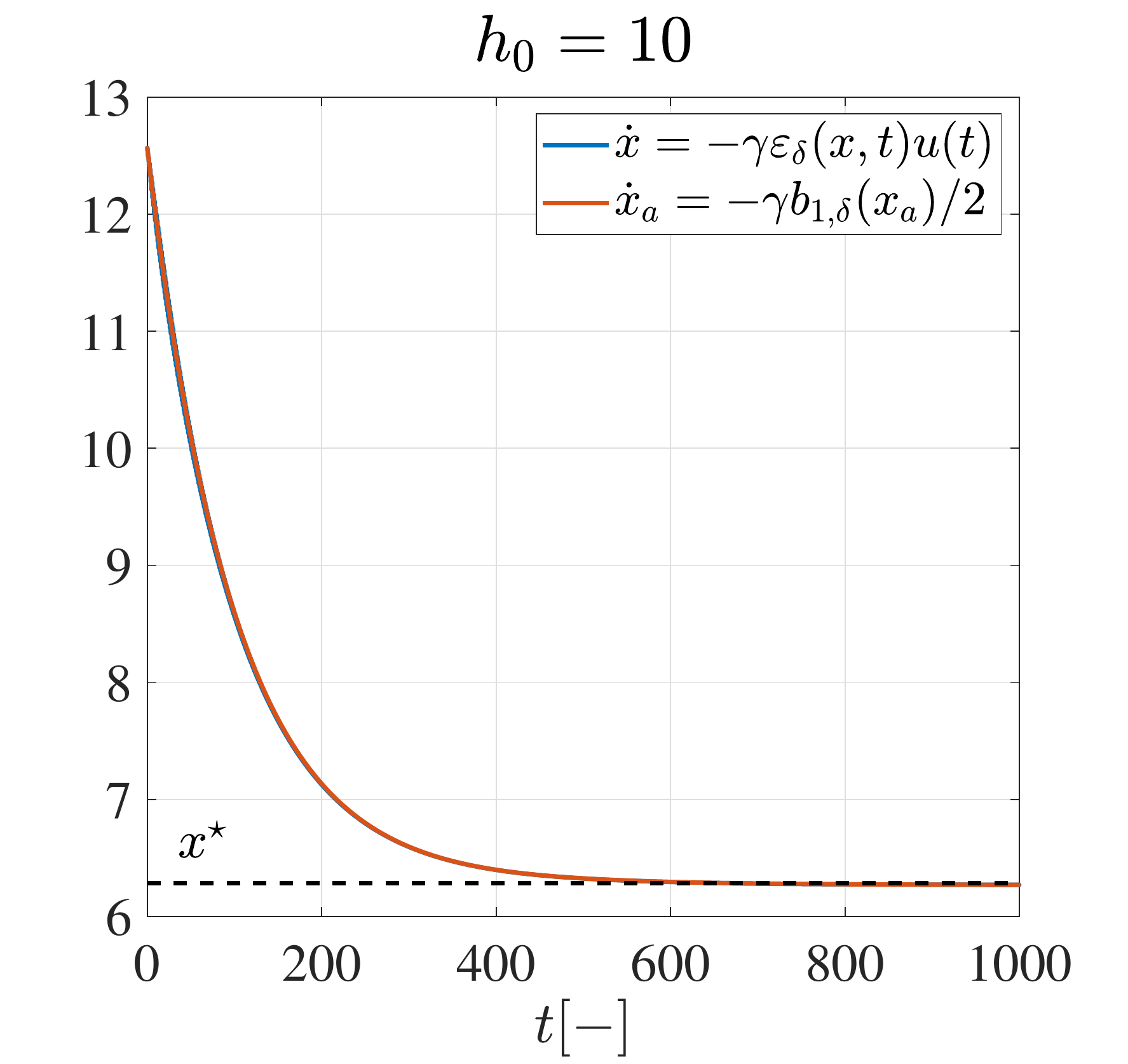}
	\includegraphics[width=0.48\columnwidth]{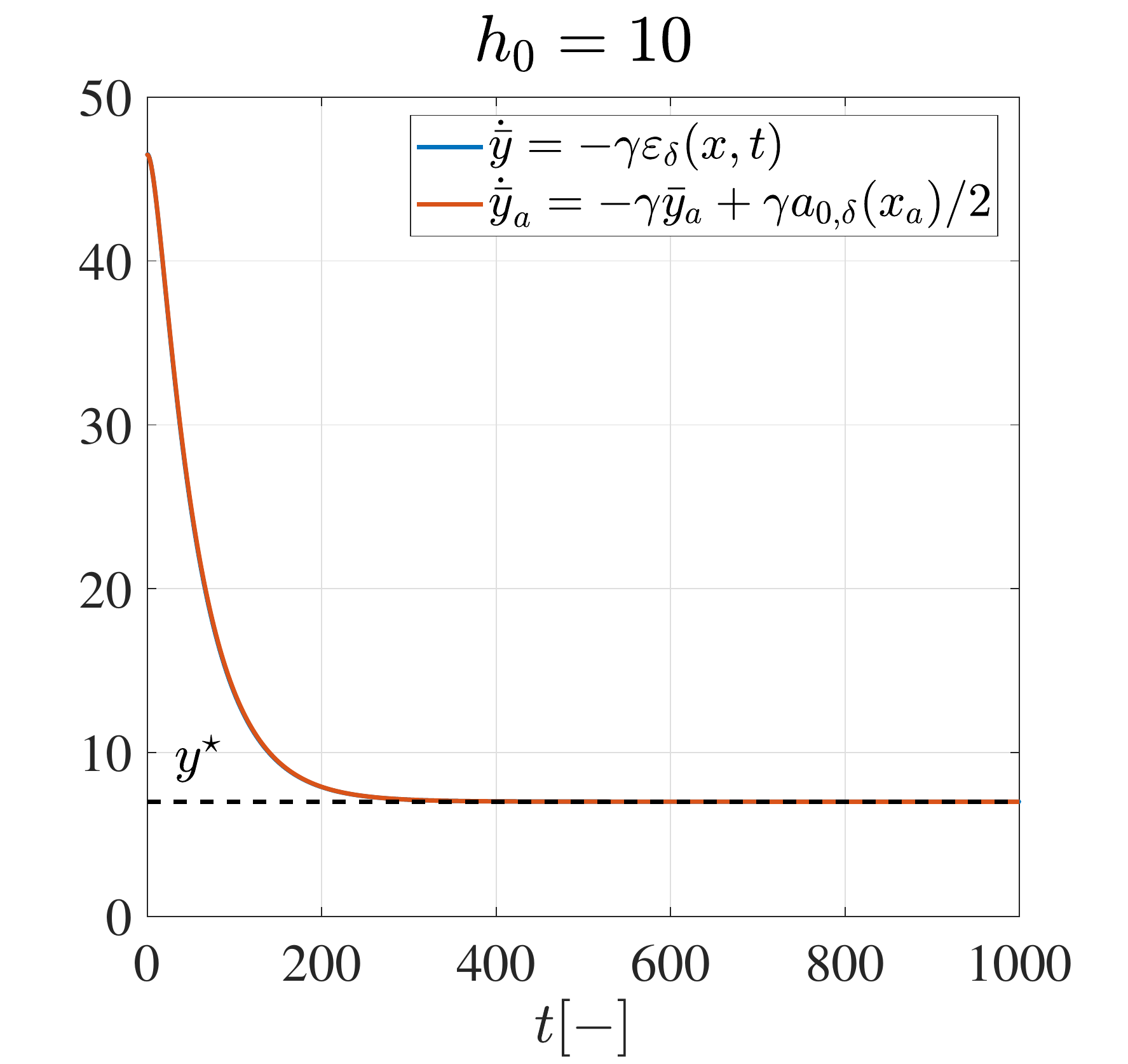}
	\includegraphics[width=0.48\columnwidth]{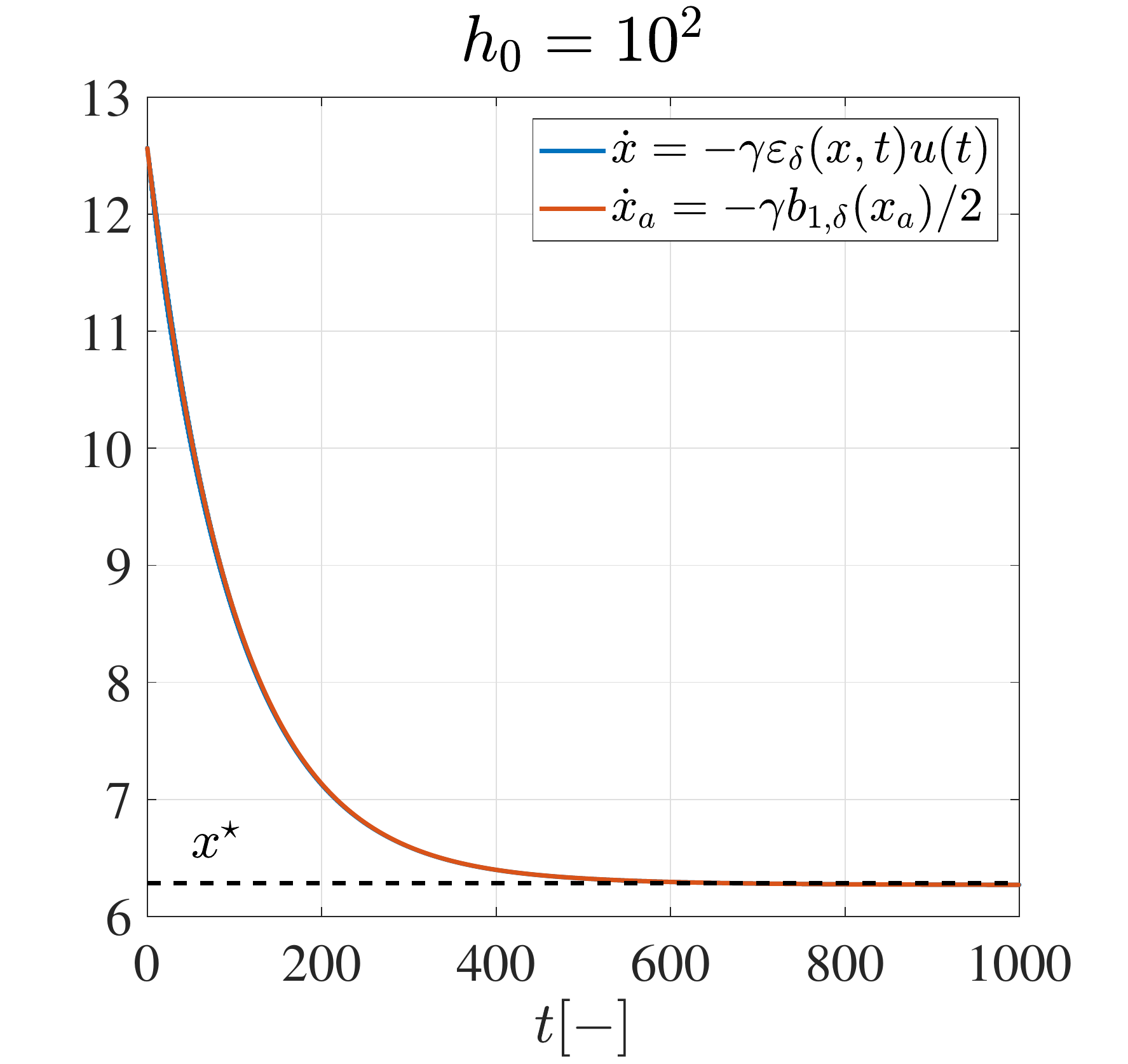}
	\includegraphics[width=0.48\columnwidth]{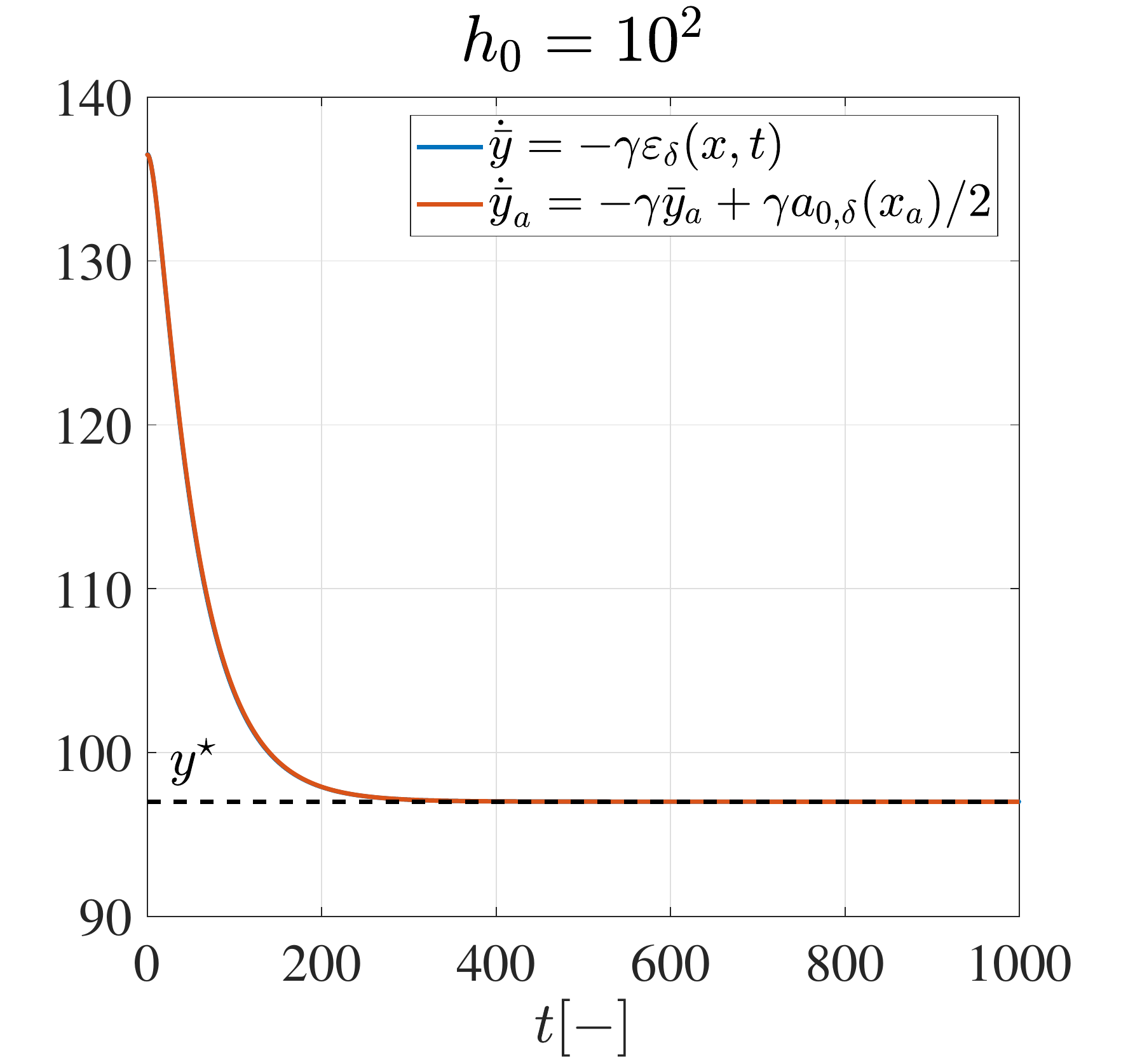}
	\includegraphics[width=0.48\columnwidth]{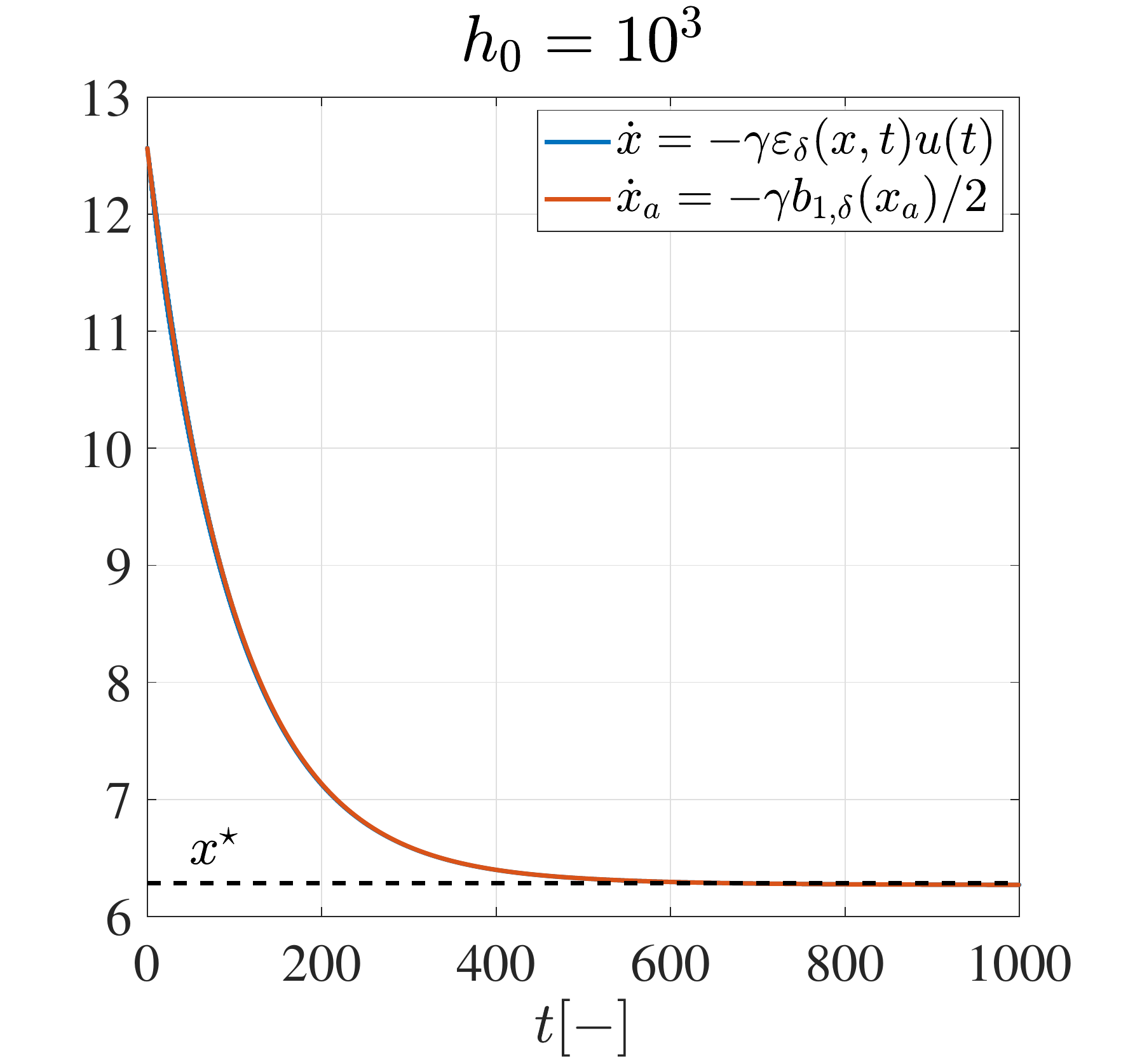}
	\includegraphics[width=0.48\columnwidth]{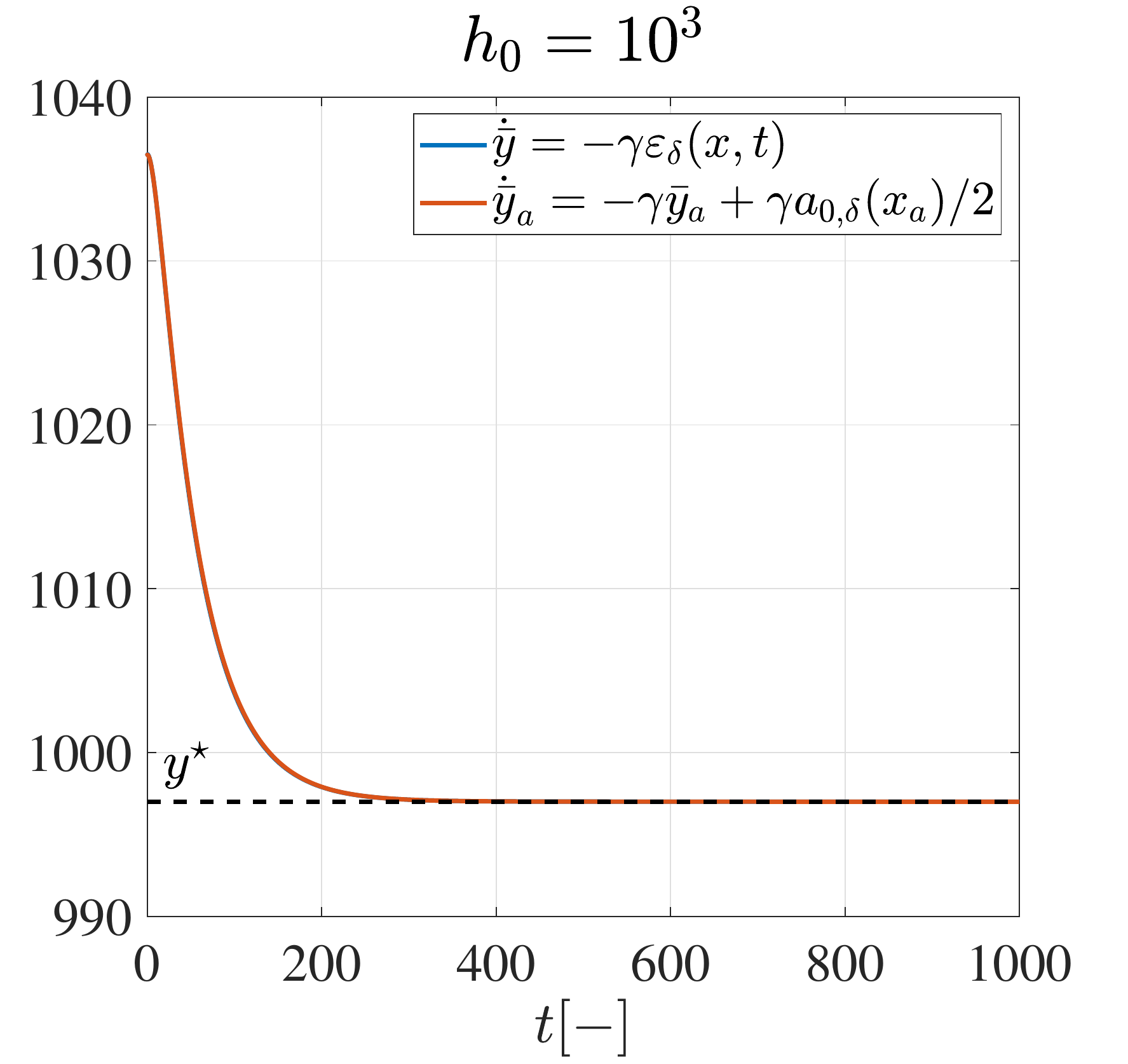}
	\caption{The behaviour of the HPF-ES does not change for increasing $M_r$. In these simulations the value of $h_0$ is increased from $-10$ to $10^3$ demonstrating that, keeping fixed $\gamma$, the HPF-ES keeps the same speed of convergence. These simulations are performed with $\gamma = \delta = 0.1$. (left) Behaviour of $x$ and $x_a$. (right) Behaviour of $\bar{y}$ and $\bar{y}_a$.}
	\label{fig:NovelES}
\end{figure}

The third group of simulations, reported in Figure \ref{fig:NovelES}, shows the performance of the HPF-ES scheme whose main feature consists of better convergence performance in case of large cost functions. 

Through the {proposed} simulations we confirm three results: a) the classic ES is able to deal with cost functions with local saddle points where, at the opposite, the average of the classic ES obtained through the Taylor expansion gets stacked (Figure \ref{fig:2ndMap_b}); b) the classic ES and its average via Taylor expansion do not converge to the same equilibrium point (Figure \ref{fig:2ndMap_c});  c) in case of large cost functions the classic ES behaves as an oscillator whose amplitude is proportional to the value of the cost function (Figure \ref{fig:Mr}(d)). Vice versa, the adoption of the high-pass filter makes the basic ES convergence rate uniform with respect to the amplitude of the cost function.  

\section{Conclusions}
\label{sec:Concl}

This paper {deals with} two well-known extremum seeking schemes {to show that they work under less restrictive assumptions}. Relying on averaging and Fourier-series arguments, it is demonstrated that these schemes are able to deal with strictly quasi-convex cost functions making the global minimiser (assumed to be unique) semi-global practically stable. Moreover, it is shown that the presence of a high-pass filter elaborating the cost function makes the tuning of the parameters independent on the cost magnitude. 
\vspace{5mm}

\bibliography{PaperBib}

\nocite{Khong2012Multidimensional}
\nocite{HAZELEGER2020109068}
\nocite{Marconi2010Robust}
\nocite{marconi2008essential}
\nocite{marconi2007output}

\end{document}